\documentclass[11pt,reqno]{amsart}

\newtheorem{proposition}{Proposition}[section]

\newtheorem{corollary}[proposition]{Corollary}
\newtheorem{theorem}[proposition]{Theorem}

\theoremstyle{definition}
\newtheorem{definition}[proposition]{Definition}

\newtheorem{remark}[proposition]{Remark}

\newcommand{\thlabel}[1]{\label{th:#1}}
\newcommand{\thref}[1]{Theorem~\ref{th:#1}}
\newcommand{\selabel}[1]{\label{se:#1}}
\newcommand{\seref}[1]{Section~\ref{se:#1}}

\newcommand{\prlabel}[1]{\label{pr:#1}}
\newcommand{\prref}[1]{Proposition~\ref{pr:#1}}
\newcommand{\colabel}[1]{\label{co:#1}}
\newcommand{\coref}[1]{Corollary~\ref{co:#1}}
\newcommand{\relabel}[1]{\label{re:#1}}
\newcommand{\reref}[1]{Remark~\ref{re:#1}}

\newcommand{\eqlabel}[1]{\label{eq:#1}}
\newcommand{\equref}[1]{(\ref{eq:#1})}

\newcommand{\Aut}{{\rm Aut}\,}

\def\NN{{\mathbb N}}

\def\ZZ{{\mathbb Z}}

\newcommand{\Cc}{\mathcal{C}}

\def\*C{{}^*\hspace*{-1pt}{\Cc}}
\def\text#1{{\rm {\rm #1}}}

\usepackage{amssymb}
\usepackage{color,amssymb,graphicx,amscd}

\begin{document}

\title[Schreier type theorems for bicrossed products] {Schreier type theorems for bicrossed products}

\author{A. L. Agore}
\address{Faculty of Engineering, Vrije Universiteit Brussel, Pleinlaan 2, B-1050 Brussels, Belgium}
\email{ana.agore@vub.ac.be and ana.agore@gmail.com}
\author{G. Militaru}
\address{Faculty of Mathematics and Computer Science, University of Bucharest, Str.
Academiei 14, RO-010014 Bucharest 1, Romania}
\email{gigel.militaru@fmi.unibuc.ro and gigel.militaru@gmail.com}
\thanks{A.L. Agore is ''Aspirant'' Fellow of the Fund for Scientific
Research–-Flanders (Belgium) (F.W.O.– Vlaanderen). G. Militaru is
supported by the CNCS - UEFISCDI grant PN-II-ID-PCE-2011-3-0039
'Hopf algebras and related topics'.}

\subjclass[2000]{20B05, 20B35, 20D06, 20D40}

\keywords{Matched pairs, bicrossed product of groups}

\begin{abstract}
We prove that the bicrossed product of two groups is a quotient of
the pushout of two semidirect products. A matched pair of groups
$(H, G, \alpha, \beta)$ is deformed using a combinatorial datum
$(\sigma, v, r)$ consisting of an automorphism $\sigma$ of $H$, a
permutation $v$ of the set $G$ and a transition map $r: G\to H$ in
order to obtain a new matched pair $\bigl(H, (G,*), \alpha',
\beta' \bigl)$ such that there exist an $\sigma$-invariant
isomorphism of groups $H\, {}_{\alpha}\!\! \bowtie_{\beta} \, G
\cong H\, {}_{\alpha'}\!\! \bowtie_{\beta'} \,(G,*)$. Moreover, if
we fix the group $H$ and the automorphism $\sigma \in \Aut(H)$
then any $\sigma$-invariant isomorphism $H\, {}_{\alpha}\!\!
\bowtie_{\beta} \, G \cong H\, {}_{\alpha'}\!\! \bowtie_{\beta'}
\, G'$ between two arbitrary bicrossed product of groups is
obtained in a unique way by the above deformation method. As
applications two Schreier type classification theorems for
bicrossed product of groups are given.
\end{abstract}

\date{}

\maketitle

\section*{Introduction}

The aim of the paper is to bring back to attention one of the most
famous open problems of group theory formulated in the first half
of the last century(\cite{Douglas}, \cite{Ore}, \cite{Redei}). It
can be seen as the dual of the more famous \textit{extension
problem} of O. L. H\"{o}lder and it is called the
\textit{factorization problem}. The statement is very simple and
tempting:

\textit{Let $H$ and $G$ be two given groups. Describe and classify
up to an isomorphism all groups $E$ that factorize through $H$ and
$G$: i.e. $E$ contains $H$ and $G$ as subgroups such that $E = H
G$ and $ H \cap G = 1$.}

Leaving aside the classification part introduced above, the first
part of the problem was formulated in 1937 by O. Ore \cite{Ore}
but it roots are much older and descend to E. Maillet's 1900 paper
\cite{Maillet}. The factorization problem has generated an
explosion of interest in group theory since the 1950's: for the
evolution and results obtained we refer to the excellent monograph
\cite{AFG}. We dare to say that the factorization problem is even
more difficult than the more popular extension problem. In the
case of two cyclic groups $H$ and $G$, not both finite, the
problem was started by L. R\'{e}dei in \cite{Redei} and finished
by P.M. Cohn in \cite{Cohn}, without the classification part
introduced above. If $H$ and $G$ are both finite cyclic groups the
problem is more difficult and seems to be still an open question,
even though J. Douglas \cite{Douglas} has devoted four papers and
over two dozen theorems to the subject. Recently, in \cite[Theorem
2.1]{ACIM} the problem was solved in the case that one of the
finite cyclic groups is of prime order. Using a theorem of
Frobenius a Schur-Zassenhaus type theorem was proven: any group
$E$ that factorizes through two finite cyclic groups, one of them
being of prime order, is isomorphic to a semidirect product
between the two cyclic groups of the same order. One of the famous
results about the factorization problem remains Ito's theorem
\cite{Ito}: let a group $E = HG$ be the product of two abelian
subgroups $H$ and $G$. Then $E$ is a metabelian group.

The converse of the factorization problem was also studied: given
a group $E$ find all \textit{exact factorizations} of it, that is,
all subgroups $H$ and $G$ of $E$ such that $E = H G$ and $ H \cap
G = 1$. Starting with the 1980's various papers dealing with this
problem were written (see \cite{Ba}, \cite{Gi}, \cite{LPS1},
\cite{LPS2}, \cite{Pr}, \cite{WW} and their list of references).
For example in \cite{LPS1} all factorizations of finite simple
groups by pairs of \textit{maximal subgroups} are described.
Derived from this problem is the following: describe and
characterize the class of (finite simple) groups that do not admit
an exact factorization between two proper subgroups. Having in
mind the abelian case such a group will be called an
\textit{indecomposable group}: the quaternion group $Q$,
$\ZZ_{p^{n}}$ for a prime integer $p$ or the alternating group
$A_{6}$ are typical examples of indecomposable groups.

An important step related to the factorization problem was the
construction of the bicrossed product $H\, {}_{\alpha}\!\!
\bowtie_{\beta} \, G$ associated to a matched pair $(H, G, \alpha,
\beta)$ given by M. Takeuchi \cite{Takeuchi}: $\alpha$ is a left
action of the group $G$ on the set $H$, $\beta$ is a right action
of the group $H$ on the set $G$ satisfying two compatibility
conditions. A group $E$ factorizes through two subgroups $H$ and
$G$ if and only if there exists a matched pair $(H, G, \alpha,
\beta)$ such that
$$
\theta : H\, {}_{\alpha}\!\! \bowtie_{\beta} \, G \rightarrow E,
\qquad \theta (h, g) = hg
$$
is an isomorphism of groups. Thus the factorization problem can be
restated in a computational manner as follows:

\textit{Let $H$ and $G$ be two given groups. Describe all matched
pairs $(H, G, \alpha, \beta)$ and classify up to an isomorphism
all bicrossed products $H\, {}_{\alpha}\!\! \bowtie_{\beta} \,
G$.}

The motivation for the above problem is triple: first of all, the
problem presents an interest in itself in group theory. On the
other hand the construction of the bicrossed product provides the
easiest way of constructing finite quantum groups \cite{masu},
hence the classification theorems from group level lead us to
classification theorems for finite quantum groups. Finally, the
bicrossed product construction at the level of groups served as a
model for similar constructions in other fields of mathematics
like: algebras \cite{cap}, coalgebras \cite{CIMZ}, groupoids
\cite{AA}, Hopf algebras \cite{Takeuchi}, locally compact groups
\cite{baaj} or locally compact quantum groups \cite{VV}, Lie
Algebras \cite{Mic} or Lie groups \cite{Kro}. Thus, the above
problem can be easily formulated for each of the above different
levels where the bicrossed product construction was made. For
instance, at the level of algebras (the bicrossed product of two
algebras is also called \emph{twisted tensor product algebra}) the
first steps were already made in the last years: the story started
with \cite[Examples 2.11]{CIMZ} where all bicrossed product
between two group algebras of dimension two are completely
described and classified. Recently, the classification of all
bicrossed product between the algebras $k^2$ and $k^m$ was
finished in  \cite{Pena} and the description of some bicrossed
products between two polynomial algebras $k[X]$ and $k[Y]$ was
started in \cite{gucci}. On the other hand, in \cite{Jara} only a
sufficient condition for the isomorphism between two bicrossed
products of algebras that fix one of the algebra is given under
the name of \textit{invariance under twisting} problem.

This paper is devoted to the classification part of the
factorization problem at the group level. Namely we shall ask the
following question: when are two bicrossed products $H\,
{}_{\alpha}\!\! \bowtie_{\beta} \, G$ and $H\, {}_{\alpha'}\!\!
\bowtie_{\beta'} \, G$ isomorphic? The organization of the paper
is the following: in \seref{1} we recall the construction of the
bicrossed product of two groups given by M. Takeuchi. It is a
generalization of the semidirect product construction for the case
when none of the factors is required to be normal. The first
natural question arises: how far is a bicrossed product from being
a semidirect product? \prref{pushout} gives the first answer to
the question: we prove that the bicrossed product of two groups is
a quotient of the pushout of two semidirect products over the
direct product of the subgroups of invariants of the actions
$\alpha$ and $\beta$. In \seref{2} we start the classification
part of the factorization problem. The main result is
\thref{deformation} : for any matched pair of groups $(H, G,
\alpha, \beta)$ and any triple $(\sigma, v, r)$, consisting of an
automorphism $\sigma$ of $H$, a permutation $v$ on the set $G$ and
a transition map $r: G\to H$ satisfying a certain compatibility
condition, a new matched pair $\bigl(H, (G,*), \alpha', \beta'
\bigl)$ is constructed such that there exists an
$\sigma$-invariant isomorphism of groups $H\, {}_{\alpha}\!\!
\bowtie_{\beta} \, G \cong H\, {}_{\alpha'}\!\! \bowtie_{\beta'}
\,(G,*)$. The importance of the result is given by the converse:
if we fix the group $H$ and the automorphism $\sigma \in \Aut(H)$,
then any $\sigma$-invariant isomorphism $H\, {}_{\alpha}\!\!
\bowtie_{\beta} \, G \cong H\, {}_{\alpha'}\!\! \bowtie_{\beta'}
\, G'$ between two arbitrary bicrossed products of groups is
obtained in a unique way by the above deformation method. As
applications in \seref{3} two Schreier type classification
theorems for bicrossed products of groups are given. They are
formulated using the language of category theory. Let $H$ and $G$
be two fixed groups: we define a category $B_{1}(H,G)$ having as
object the set of all matched pairs $(H, G, \alpha, \beta)$ and
morphisms are defined as morphisms between two bicrossed products
that fix one of the groups. \thref{sch1} gives a bijection between
the set of objects of the skeleton of the category $B_{1}(H,G)$
and a certain pointed set $K^{2}(H, G)$ that will play for the
classification problem of bicrossed products the same role as the
second cohomology group does for the classification of the
extension problem. Returning to the question of how far a
bicrossed product is from being a semidirect product,
\coref{lake1} and \coref{lake2} give two necessary and sufficient
conditions for a bicrossed product to be isomorphic to a
semidirect product of groups in the category $B_{1}(H,G)$.
\thref{sch22} is the second Schreier type theorem for bicrossed
products: this time we fix two groups $H$, $G$ and $\beta: G
\times H \rightarrow G$ a right action of the group $H$ on the set
$G$ and the classification theorem is more restrictive than the
one given in \thref{sch1}. In the last section we give some
examples: we compute and count explicitly the set of all matched
pairs $(C_3, C_m, \alpha, \beta)$, where $C_m$ is a cyclic group
of order $m$, and the pointed set $K^{2} (C_{3}, C_{6})$
constructed in \thref{sch1} is shown to have three elements.

\section{Preliminaries}\selabel{1}
Let us fix the notation that will be used throughout the paper.
Let $H$ and $G$ be two groups and $\alpha : G \times H \rightarrow
H$ and $\beta : G \times H \rightarrow G$ two maps. We use the
notation
$$\alpha (g, h) = g\triangleright h \quad \text{and} \quad
\beta (g, h) = g\triangleleft h$$ for all $g\in G$ and $h\in H$.
The map $\alpha$ (resp. $\beta$) is called trivial if
$g\triangleright h = h$ (resp. $g\triangleleft h = g$) for all
$g\in G$ and $h\in H$. We recall that $\alpha$ is an action as
automorphism if it is a left action of the group $G$ on the set
$H$ and $g \triangleright (h_1 h_2) = (g \triangleright h_1) (g
\triangleright h_2)$, for all $g\in G$, $h_1$, $h_2 \in H$.
Similarly, $\beta$ is an action as automorphism if it is a right
action of the group $H$ on the set $G$ and $(g_1 g_2) \lhd h =
(g_1 \lhd h) (g_2 \lhd h)$ for all $g_1$, $g_2 \in G$ and $h\in
H$. $\Aut(H)$ is the group of automorphisms of $H$ and $C_n$ is
the cyclic group of order $n$.

Let $H$ and $G$ be two groups with the multiplications $m_{H}: H
\times H\rightarrow H$, $m_{G}:G \times G\rightarrow G$, units
$1_H$ and respectively $1_G$ and $R: G \times H \rightarrow H
\times G$ a map. We shall define a new multiplication on the set
$H \times G$ using $R$ instead of the usual flip $\tau: G\times H
\to H\times G$, $\tau (g, h) = (h, g)$ as follows:
$$m_{H \times G, R}: H \times G \times H \times G \rightarrow H
\times G, \qquad m_{H \times G, R}:= (m_{H} \times m_{G})\circ (I
\times R \times I)$$

Let $\alpha := \pi_{1} \circ R : G \times H \rightarrow H$, $\beta
:=\pi_{2} \circ R: G \times H \rightarrow G$, where $\pi_{i}$ is
the projection on the $i$-component; we shall denote $\alpha (g,
h) = g \rhd h$ and $\beta (g, h) = g \lhd h$, for all $g\in G$ and
$h\in H$. Then $R (g, h) = (g \rhd h, g \lhd h)$ and the
multiplication $m_{H \times G, R}$ on $H \times G$ can be
explicitly written as follows:
\begin{equation}
(h_{1}, g_{1}) \cdot_{R} (h_{2}, g_{2}) = \bigl( h_{1}(g_{1} \rhd
h_{2}), \, (g_{1} \lhd h_{2})g_{2}\bigl)
\end{equation}
for all $h_{1}$, $h_{2} \in H$ and $g_{1}$, $g_{2} \in G$.

It can be easily shown that $(H \times G, m_{H \times G, R})$ is a
group with $(1_{H}, 1_{G})$ as a unit if and only if $(H, G,
\alpha, \beta)$ is a \textit{matched pair} in the sense of
Takeuchi (\cite{Takeuchi}): i.e. $\alpha$ is a left action of the
group $G$ on the set $H$, $\beta$ is a right action of the group
$H$ on the set $G$ and the following two compatibility conditions
hold:
\begin{equation}\eqlabel{2}
g \rhd (h_{1} h_{2}) = (g \rhd h_{1})\bigl((g \lhd h_{1}) \rhd
h_{2}\bigl)
\end{equation}
\begin{equation}\eqlabel{3}
(g_{1} g_{2}) \lhd h = \bigl(g_{1} \lhd (g_{2} \rhd h)\bigl)(g_{2}
\lhd h)
\end{equation}
for all $h$, $h_{1}$, $h_{2} \in H$ and $g$, $g_{1}$, $g_{2} \in
G$. It follows from \equref{2} and \equref{3} that:
\begin{equation}\eqlabel{4}
g \rhd 1_H = 1_H \quad {\rm and}\quad 1_G \lhd h = 1_G
\end{equation}
for all $h \in H$ and $g \in G$.

If $(H, G, \alpha, \beta)$ is a matched pair, the new group
obtained on the set $H \times G$ will be denoted by $H\,
{}_{\alpha}\!\! \bowtie_{\beta} \, G = H\bowtie \, G$ and will be
 called the \textit{bicrossed product} ({\sl knit product} or
{\sl Zappa-Sz\' ep product}) of $H$ and $G$. We note that $i_H :
H\to H\bowtie \, G$, $i_H (h) = (h, 1)$ and $i_G : G\to H\bowtie
\, G$, $i_G (g) = (1, g)$, for all $h\in H$, $g\in G$ are
morphisms of groups and hence $H \times \{1\}\cong H$ and $\{1\}
\times G \cong G$ are subgroups of $H\bowtie \, G$. Moreover,
every element $(h, g)$ of $H\bowtie \, G$ can be written uniquely
as a product of an element of $H \times \{1\}$ and of an element
of $\{1\} \times G$ as follows:
\begin{equation}
(h,g) = (h,1) \cdot (1,g)
\end{equation}
Conversely, this observation characterizes the bicrossed product.
Let $E$ be a group $H$, $G\leq E$ be subgroups such that any
element of $E$  can be written uniquely as a product of an element
of $E$ and an element of $G$. Then there exists a matched pair
$(H, G, \alpha, \beta)$ such that
$$
\theta : H\bowtie \, G \ \rightarrow E, \qquad \theta (h, g) = hg
$$
is a group isomorphism (\cite{Takeuchi}). The maps $\alpha$ and
$\beta$ play a symmetric role: if $(H, G, \alpha, \beta)$ is a
matched pair then we can construct a new matched pair $(G, H,
\tilde{\alpha}, \tilde{\beta} )$ such that there exists a
canonical isomorphism of groups  $H\, _{\alpha}\bowtie_{\beta} \,
G  \cong G\, _{\tilde{\alpha}}\bowtie_{\tilde{\beta}} \, H$
(\cite[Proposition 2.5]{ACIM}).

\begin{remark}\relabel{2.4.90}
Let $H$ and $G$ be two groups and $\beta : G \times H \rightarrow
G$ the trivial action. Then $(H, G, \alpha, \beta)$ is a matched
pair if and only if $\alpha : G \times H \rightarrow H$ is an
action of $G$ on $H$ as group automorphisms. In this case the
bicrossed product $H\bowtie \, G$ is exactly the left version of
the semidirect product $H {}_{\alpha}\ltimes G$.

Assume now that the map $\alpha$ is the trivial action. Then $(H,
G, \alpha, \beta)$ is a matched pair if and only if $\beta$ is a
right action of $H$ on $G$ as group automorphisms. Is this case
the bicrossed product $H\bowtie \, G$ is exactly the right version
of the semidirect product $H\rtimes_{\beta} \, G$. It can be
easily proved that $H\rtimes_{\beta} \, G \cong G
{}_{\varphi}\ltimes H$, where $\varphi = \varphi_{\beta}$ is the
action of $H$ on $G$ as group automorphisms given by
$$
\varphi : H \rightarrow \Aut (G), \qquad \varphi (h)(g) = \Bigl(
g^{-1} \triangleleft {h^{-1}} \Bigl)^{-1}
$$
for all $h\in H$ and $g\in G$ (\cite[Remark 2.6]{ACIM}).

A matched pair $(H, G, \alpha, \beta)$ is called \textit{proper}
if $\alpha$ and $\beta$ are both nontrivial actions.
\end{remark}

The above Remark shows that the semidirect product is a special
case of the bicrossed product construction. It is therefore
natural to ask the converse: \textit{Can a bicrossed product be
obtained from semidirect products of groups?} In what follows we
shall give a first answer to this question: a bicrossed product
can be obtained as a quotient of a pushout of two semidirect
products of groups.

Let $(H, G, \alpha, \beta)$ be a matched pair and let us denote by
${\rm Fix}(H)$ and ${\rm Fix}(G)$ the invariants of the two
actions $\alpha$ and $\beta$ :
$$
{\rm Fix} (H) : = \{h \in H \mid g \rhd h = h, \, \forall g \in
G\}, \quad {\rm Fix}(G) : = \{g \in G \mid g \lhd h = g, \,
\forall h \in H\}
$$
Using the compatibility conditions \equref{2} and \equref{3} we
shall prove that ${\rm Fix}(H)$ is a subgroup of $H$ and ${\rm
Fix}(G)$ a subgroup of $G$. Indeed, from \equref{4} we obtain that
$1_H \in {\rm Fix}(H)$ and for $h_{1}$, $h_{2} \in {\rm Fix}(H)$
we have:
$$ g \rhd (h_{1}h_{2}) \stackrel{\equref{2}} {=}
(g \rhd h_{1})\bigl((g \lhd h_{1}) \rhd h_{2}\bigl) = h_{1}h_{2}
$$
i.e. $h_{1}h_{2} \in {\rm Fix}(H)$. On the other hand:
$$1_H
\stackrel{\equref{4}} {=} g \rhd 1_H = g \rhd (h_{1}^{-1}h_{1})
\stackrel{\equref{2}} {=} (g \rhd h_{1}^{-1})\bigl((g \lhd
h_{1}^{-1}) \rhd h_{1}\bigl) = (g \rhd h_{1}^{-1})h_{1}$$

Thus $g \rhd h_{1}^{-1} = h_{1}^{-1}$, i.e. $h_{1}^{-1} \in {\rm
Fix}(H)$. In a similar way we can show that ${\rm Fix}(G)$ is a
subgroup of $G$. Using the compatibility condition \equref{2} we
obtain that the map given by:
$$\varphi_{\rhd}: {\rm Fix}(G) \rightarrow {\rm Aut}(H), \quad
\varphi_{\rhd}(g)(h) := g \rhd h$$ for all $g \in {\rm Fix}(G)$,
$h\in H$ is a morphism of groups. Thus we can construct the left
version of the semidirect product associated to the triple $(H, \,
{\rm Fix}(G), \, \varphi_{\rhd})$: that is $H
{}_{\varphi_{\rhd}}\!\! \ltimes {\rm Fix}(G) : = H\times {\rm
Fix}(G)$ with the multiplication:
$$(h, g) (h', g') = \bigl(h(g \rhd h'), \, gg'\bigl)$$
for all $h$, $h' \in H$ and $g$, $g' \in {\rm Fix}(G)$. Similarly,
using \equref{3} we obtain that the map given by:
$$\psi_{\lhd}: {\rm Fix}(H) \rightarrow {\rm Aut}(G), \quad \psi_{\lhd}(h)(g):= g \lhd h$$
for all $h \in {\rm Fix}(H)$, $g\in G$ is a morphism of groups and
we can construct the right version of the semidirect product
associated to the triple $(G, \, {\rm Fix}(H), \, \psi_{\lhd})$:
i.e. ${\rm Fix}(H) \rtimes_{\psi_{\lhd}} G := {\rm Fix}(H) \times
G$ with the multiplication:
$$(h, g)(h', g') = \bigl(hh', \, (g \lhd h')g'\bigl)$$
for all $h$, $h' \in {\rm Fix}(H)$ and $g$, $g' \in G$. Moreover,
the inclusion maps
$$\overline{i}: {\rm Fix}(H) \times {\rm Fix}(G) \hookrightarrow H
{}_{\varphi_{\rhd}}\!\! \ltimes {\rm Fix}(G) \quad {\rm and} \quad
\overline{j} : {\rm Fix}(H) \times {\rm Fix}(G) \hookrightarrow
{\rm Fix}(H) \rtimes_{\psi_{\lhd}} G$$ are morphisms of groups by
straightforward verifications.

On the other hand we can easily prove that the canonical
inclusions
$$i: H {}_{\varphi_{\rhd}}\!\! \ltimes {\rm Fix}(G)
\hookrightarrow H\, {}_{\alpha}\!\! \bowtie_{\beta} \, G, \quad
i(h,g) = (h,g)$$ and
$$j: {\rm Fix}(H) \rtimes_{\psi_{\lhd}} G
\hookrightarrow H\, {}_{\alpha}\!\! \bowtie_{\beta} \, G, \quad
j(h,g) = (h,g)$$ are morphisms of groups. Indeed for $h$, $h'\in
H$ and $g$, $g' \in {\rm Fix}(G)$ we have:
$$ i(h, g) \cdot i(h', g') = \bigl( h (g \rhd h'), (g \lhd h')
g'\bigl)\, \stackrel{g \in {\rm Fix}(G)} {=} \bigl( h(g \rhd h'),
gg'\bigl) = i\bigl((h, g) (h',g')\bigl)$$

Thus the two semidirect products constructed above, $H
{}_{\varphi_{\rhd}}\!\! \ltimes {\rm Fix}(G)$ and ${\rm Fix}(H)
\rtimes_{\psi_{\lhd}} G$, are subgroups of the bicrossed product
$H\, {}_{\alpha}\!\! \bowtie_{\beta} \, G$. To conclude, we
obtained a commutative diagram in the category of groups
\begin{equation}\eqlabel{sfibr}
\begin {CD}
Fix(H)\times Fix(G) @>\overline{j}>> Fix(H)\rtimes_{\psi} G\\
@VV\overline{i}V @VVjV\\
H {}_\varphi \ltimes Fix(G)@>i>> H {}_\alpha \bowtie_\beta G
\end{CD}
\end{equation}

Using the construction of the pullback in the category of groups
it follows that the pair $({\rm Fix}(H) \times {\rm Fix}(G),
(\overline{i}, \overline{j}))$ is a pullback of the morphisms $i:
H {}_{\varphi_{\rhd}}\!\! \ltimes {\rm Fix}(G) \hookrightarrow H\,
{}_{\alpha}\!\! \bowtie_{\beta} \, G$ and $j: {\rm Fix}(H)
\rtimes_{\psi_{\lhd}} G \hookrightarrow H\, {}_{\alpha}\!\!
\bowtie_{\beta} \, G$.

\begin{proposition}\prlabel{pushout}
Let $(H, G, \alpha, \beta)$ be a matched pair of groups and
$\bigl(X, (\varphi, \psi)\bigl)$ be the pushout in the category of
groups of the diagram
$$
\begin {CD}
Fix(H)\times Fix(G) @>\overline{j}>> Fix(H)\rtimes_{\psi} G\\
@VV\overline{i}V @VV\varphi V\\
H {}_\varphi \ltimes Fix(G)@>\psi >> X
\end{CD}
$$
Then the bicrossed product $H\, {}_{\alpha}\!\! \bowtie_{\beta} \,
G$ is isomorphic to a quotient group of $X$.
\end{proposition}

\begin{proof}
The diagram \equref{sfibr} is commutative and $\bigl(X, (\varphi,
\psi)\bigl)$ is the pushout of the pair $(\overline{i}, \,
\overline{j})$: thus there exists an unique morphism of groups
$\theta : X \rightarrow H\, {}_{\alpha}\!\! \bowtie_{\beta} \, G$
such that $\theta \circ \psi = i$ and $\theta \circ  \varphi = j$.
Let $(h, g) \in H\, {}_{\alpha}\!\! \bowtie_{\beta} \, G$: as $(h,
1_G) \in H {}_{\varphi_{\rhd}}\!\! \ltimes Fix(G)$ and $(1_H, g)
\in Fix(H) \rtimes_{\psi_{\lhd}} G$ we obtain
$$(h, g) = (h, 1_G)(1_H, g) = i(h, 1_G) j(1_H, g) =
\theta \bigl( \psi(h, 1_G)\bigl) \theta \bigl(\varphi(1_H, g)
\bigl) = \theta \bigl(\psi(h, 1_G) \varphi(1_H, g) \bigl)$$ that
is $\theta$ is surjective. Thus $H\, {}_{\alpha}\!\!
\bowtie_{\beta} \, G$ is a quotient group of $X$.
\end{proof}

We end the section with a problem that can be of interest for a
further study:

\textit{Let "P" be a property in the category of groups. Give a
necessary and sufficient condition such that $H\, {}_{\alpha}\!\!
\bowtie_{\beta} \, G$ has the property "P".}

In the following we give an example in the case that "P" is the
property of being abelian or cyclic.

\begin{proposition} Let $(H, G, \alpha, \beta)$ be a matched pair of
groups. Then:
\begin{enumerate}
\item The center of the bicrossed product $H\, {}_{\alpha}\!\!
\bowtie_{\beta} \, G$ is given by:
$$Z\bigl(H\, {}_{\alpha}\!\! \bowtie_{\beta} \, G \bigl) = \{(h, g)
\in {\rm Fix}(H)\times {\rm Fix}(G) \mid g \rhd x = h^{-1}xh, \, y
\lhd h = gyg^{-1}, \forall x \in H, y \in G \}$$
 \item $H\, {}_{\alpha}\!\! \bowtie_{\beta} \, G$ is an abelian
group if and only if $H$ and $G$ are abelian groups and $\alpha$
and $\beta$ are the trivial actions;

\item $H\,{}_{\alpha}\!\!\bowtie_{\beta} \, G$ is a cyclic group
if and only if $\alpha$ and $\beta$ are the trivial actions and
$H$, $G$ are finite cyclic groups of coprime orders.
\end{enumerate}
\end{proposition}

\begin{proof} An element $(h, g) \in H\, {}_{\alpha}\!\!
\bowtie_{\beta} \, G$ belongs to the center of the group if and
only if $(h, g) (x, 1) = (x, 1) (h, g)$ and $(h, g) (1, y) = (1,
y) (h, g)$, for all $x \in H$ and $y \in G$. This is equivalent to
$h(g \rhd x) = xh$, $g \lhd x = g$, $y \rhd h = h$ and $(y \lhd
h)g = gy$, for all $x \in H$, $y \in G$. Hence $h \in {\rm
Fix}(H)$, $g \in {\rm Fix}(G)$, $g \rhd x = h^{-1}xh$, $y \lhd h =
gyg^{-1}$ for all $x \in H$, $y \in G$. (2) follows from (1) and
(3) follows from (2) and the Chinese lemma: a direct product of
two groups is a cyclic group if and only if they are finite,
cyclic of coprime order.
\end{proof}

\section{Deformation of a matched pair}\selabel{2}

Let $H$ be a group and $\sigma \in \Aut(H)$ an automorphism of
$H$. We define the category $\mathcal{C}(H, \sigma)$ as follows:
an object of $\mathcal{C}(H, \sigma)$ is a triple $(G, \alpha,
\beta)$ such that $(H, G, \alpha, \beta)$ is a matched pair of
groups. A morphism $\psi : (G', \alpha', \beta') \rightarrow (G,
\alpha, \beta)$ in $\mathcal{C}(H, \sigma)$ is a morphism of
groups $\psi: H\, {}_{\alpha'}\!\! \bowtie_{\beta'} \, G'
\rightarrow H\, {}_{\alpha}\!\! \bowtie_{\beta} \, G$ such that
the following diagram
\begin{equation}\eqlabel{D1}
\begin{CD}
H@>i_H>> H\, {}_{\alpha'}\bowtie_{\beta'}G' \\
@VV\sigma V @VV\psi V\\
H@>i_H>> H _{\alpha} \bowtie_{\beta} G
\end{CD}
\end{equation}
is commutative. A (iso)morphism $\psi: H\, {}_{\alpha'}\!\!
\bowtie_{\beta'} \, G' \rightarrow H\, {}_{\alpha}\!\!
\bowtie_{\beta} \, G$ in the category $\mathcal{C}(H, \sigma)$
will be called a \textit{$\sigma$-invariant (iso)morphism} between
the two bicrossed products.

The following key proposition describes explicitly the morphisms
of $\mathcal{C}(H, \sigma)$ and gives a necessary and sufficient
condition for two bicrossed products $H\, {}_{\alpha'}\!\!
\bowtie_{\beta'} \, G'$ and $H\, {}_{\alpha}\!\! \bowtie_{\beta}
\, G$ to be isomorphic in the category $\mathcal{C}(H, \sigma)$.
If $G'$ is a new group we shall denote by "$*$" the multiplication
of $G'$ and $\alpha'(g',h) = g' \rhd' h$, $\beta' (g',h) = g'
\lhd' h$, for all $g'\in G'$ and $h\in H$.

\begin{proposition}\prlabel{1}
Let $H$ be a group, $\sigma \in \Aut(H)$ and $(H, G, \alpha,
\beta)$, $(H, G', \alpha', \beta')$ two matched pairs. There
exists a one to one correspondence between the set of all
morphisms $\psi: H\, {}_{\alpha'}\!\! \bowtie_{\beta'} \, G'
\rightarrow H\, {}_{\alpha}\!\! \bowtie_{\beta} \, G$ in the
category $\mathcal{C}(H, \sigma)$ and the set of all pairs
$(r,v)$, where $r: G' \rightarrow H$, $v: G' \rightarrow G$ are
two maps such that:
\begin{eqnarray}
\sigma(g' \rhd' h)r(g' \lhd' h) &=& r(g')\bigl(v(g') \rhd
\sigma(h)\bigl)\eqlabel{p1} \\
v(g' \lhd' h) &=& v(g') \lhd \sigma(h) \eqlabel{p2} \\
r(g_{1}' * g_{2}') &=& r(g_{1}') \bigl( v(g_{1}') \rhd
r(g_{2}') \bigl) \eqlabel{p3} \\
v(g_{1}' * g_{2}') &=& \bigl( v(g_{1}') \lhd r(g_{2}')\bigl)
v(g_{2}') \eqlabel{p4}
\end{eqnarray}
for all $g', g_{1}', g_{2}' \in G'$, $h \in H$. Through the above
bijection $\psi$ is given by
\begin{equation}\eqlabel{p5}
\psi(h, \, g') = \bigl(\sigma(h)r(g'), \, v(g')\bigl)
\end{equation}
for all $h \in H$, $g' \in G'$. Moreover, $\psi : H\,
{}_{\alpha'}\!\! \bowtie_{\beta'} \, G' \rightarrow H\,
{}_{\alpha}\!\! \bowtie_{\beta} \, G$ is an isomorphism in
$\mathcal{\mathcal{C}}(H,\sigma)$ if and only if the map $v: G'
\to G$ is bijective.
\end{proposition}

\begin{proof}
A morphism of groups $\psi : H\, {}_{\alpha'}\!\! \bowtie_{\beta'}
\, G' \rightarrow H\, {}_{\alpha}\!\! \bowtie_{\beta} \, G$ that
makes the diagram \equref{D1} commutative is uniquely defined by
two maps $r = r_{\psi}: G' \rightarrow H$, $v = v_{\psi}:G'
\rightarrow G$ such that $\psi(1,g') = \bigl(r(g'),v(g')\bigl)$
for all $g' \in G'$. In this case $\psi$ is given by :
$$
\psi(h, g') = \psi\bigl((h,1)\cdot (1,g')\bigl) = \bigl(\sigma(h),
1\bigl) \cdot \bigl(r(g'), v(g')\bigl) = \bigl(\sigma(h)r(g'),
v(g')\bigl)
$$
for all $h \in H$ and $g' \in G'$. As $\psi(1,1)=(1,1)$ we obtain
that $r(1)=1$ and $v(1)=1$.

We shall prove now that $\psi$ is a morphism of groups if and only
if the compatibility conditions \equref{p1} - \equref{p4} hold for
the pair $(r,v)$. It is enough to check the condition $\psi(xy) =
\psi(x) \psi(y)$ only for generators $x$, $y \in \bigl(H\times
\{1\}\bigl)\cup \bigl(\{1\}\times G' \bigl)$ of the bicrossed
product $H\, {}_{\alpha'}\!\! \bowtie_{\beta'} \, G'$. Since
$\sigma$ is an automorphism of $H$, we have to check only for $x =
(1, g')$, $y = (h, 1)$ and $x = (1, g_{1}')$, $y = (1, g_{2}')$.
The condition $\psi\bigl((1,g')(h,1)\bigl) = \psi(1,g')\psi(h,1)$
is equivalent to \equref{p1} - \equref{p2} and the condition
$\psi\bigl((1,g_{1}')(1,g_{2}')\bigl) =
\psi(1,g_{1}')\psi(1,g_{2}')$ is equivalent to \equref{p3} -
\equref{p4}. Note that the normalization conditions $v(1) = 1$ and
$r(1) = 1$ where used to obtain \equref{p1} and \equref{p4}.

Conversely, the normalization conditions follow from \equref{p1} -
\equref{p4} in the following manner: first, for $g' = 1$ in
\equref{p1} we obtain $\sigma(h)r(1) = r(1)\bigl(v(1) \rhd
\sigma(h)\bigl) $ for all $h \in H$. Since $\sigma$ is an
automorphism we have :
\begin{equation}\eqlabel{th2}
hr(1) = r(1)\bigl(v(1) \rhd h\bigl)
\end{equation}
for all $h \in H$. Now let $g_{1}'=1$ in \equref{p3} to obtain:
\begin{eqnarray*}
r(g_{2}') = r(1)\bigl(v(1) \rhd r(g_{2}')\bigl)
\stackrel{\equref{th2}} {=} r(g_{2}')r(1)
\end{eqnarray*}
thus $r(1) = 1$. Finally we let $g_{2}' = 1$ in \equref{p4} to
obtain $v(g_{1}') = v(g_{1}')v(1)$, thus $v(1) = 1$.

It remains to be proven that $\psi$ given by \equref{p5} is an
isomorphism if and only if $v: G'\to G$ is a bijective map. Assume
first that $\psi$ is an isomorphism. Then $v$ is surjective and
for $g_{1}'$, $g_{2}' \in G'$ such that $v(g_{1}') = v(g_{2}')$ we
have:
$$
\psi(1, g_{2}') = \bigl(r(g_{2}'), v(g_{2}')\bigl) =
\bigl(r(g_{2}'), v(g_{1}')\bigl) =
\psi\bigl(\sigma^{-1}(r(g_{2}'))\sigma^{-1}(r(g_{1}')^{-1}),
g_{1}'\bigl)
$$
Hence $g_{1}' = g_{2}'$ and $v$ is injective. Conversely, assume
that $v$ is bijective. If $\psi(h, g') = (1, 1)$ we obtain that
$\sigma(h)r(g') = 1$ and $v(g') = 1 = v(1)$. It follows from here
that $g' = 1$ and $\sigma(h) = 1 = \sigma(1)$, i.e. $h = 1$. Hence
$\psi$ is injective. Let $(h, g) \in H\, {}_{\alpha}\!\!
\bowtie_{\beta} \, G$ and $g' \in G'$ such that $v(g') = g$. Then
$\psi\bigl(\sigma^{-1}(h)\sigma^{-1}(r(g')^{-1}), g'\bigl) = (h,
g)$ i.e. $\psi$ is an isomorphism of groups.
\end{proof}

We shall prove now the main result of this section:

\begin{theorem}\textbf{(Deformation of a matched pair)}\thlabel{deformation}
Let $(H, G, \alpha, \beta)$ be a matched pair of groups, $(\sigma,
v, r)$ be a triple where $\sigma \in \Aut(H)$, $v: G \to G$ is a
bijective map, $r: G \rightarrow H$ is a map such that $v (1_G) =
1_G$, $r (1_G) = 1_H$ and the following compatibility condition:
\begin{equation}\eqlabel{crossed}
r \circ v^{-1}\bigl((v(g_{1}) \lhd r(g_{2}))v(g_{2})\bigl) =
r(g_{1})\bigl(v(g_{1}) \rhd r(g_{2})\bigl)
\end{equation}
holds for all $g_{1}, g_{2} \in G$. On the set $G$ we define a new
multiplication $*$ and two new actions $\beta' : G\times H
\rightarrow G$, $\alpha' : G\times H \rightarrow H$ given by:
\begin{eqnarray}
g_{1} * g_{2} &:=& v^{-1}\bigl((v(g_{1}) \lhd
r(g_{2}))v(g_{2})\bigl) \eqlabel{def1}\\
g \lhd' h &:=& v^{-1}\bigl(v(g) \lhd \sigma(h)\bigl)
\eqlabel{def2}\\
g \rhd' h &:=& \sigma^{-1}(r(g))\sigma^{-1}(v(g) \rhd \sigma(h))
\sigma^{-1}\bigl(r \circ v^{-1}(v(g) \lhd \sigma(h))^{-1}\bigl)
\eqlabel{def3}
\end{eqnarray}
for all $g_1$, $g_2$, $g \in G$ and $h\in H$. Then:
\begin{enumerate}
\item $(G,*)$ is a group structure on the set $G$ with $1_{G}$ as
a unit;
\item $\bigl(H, (G,*), \alpha', \beta'\bigl)$ is a matched
pair of groups and
$$\psi : H\, {}_{\alpha'}\!\! \bowtie_{\beta'} \, (G,*)
\rightarrow H\, {}_{\alpha}\!\! \bowtie_{\beta} \, G, \quad
\psi(h,g) := \bigl(\sigma(h)r(g), v(g)\bigl)$$ is a
$\sigma$-invariant isomorphism of groups. \item Any
$\sigma$-invariant isomorphism of groups $H\, {}_{\alpha'}\!\!
\bowtie_{\beta'} \, G \cong H\, {}_{\alpha}\!\! \bowtie_{\beta} \,
G$ arises as above.
\end{enumerate}
\end{theorem}

\begin{proof}
(1) Let $g \in G$. Then
$$g * 1_G = v^{-1}\bigl((v(g) \lhd
r(1_G))v(1)\bigl) = v^{-1}\bigl(v(g)\bigl) = g$$ and
$$1_G * g =
v^{-1}\bigl((v(1_G) \lhd r(g))v(g)\bigl) = v^{-1}\bigl(v(g)\bigl)
= g$$ Hence $1_{G}$ is a unit for $*$. Let $g_{1}$, $g_{2}$,
$g_{3} \in G$. Then:
\begin{eqnarray*}
v \bigl(\underline{(g_{1} * g_{2})} * g_{3}\bigl)
&\stackrel{\equref{def1}} {=}& v\bigl[\underline{
v^{-1}\bigl((v(g_{1}) \lhd r(g_{2}))v(g_{2})\bigl) * g_{3} }\bigl]\\
&\stackrel{\equref{def1}} {=}& \bigl[ \underline{\bigl((v(g_{1})
\lhd
r(g_{2}))v(g_{2}) \bigl) \lhd r(g_{3})}\bigl]v(g_{3})\\
&\stackrel{\equref{3}} {=}& \bigl[\bigl(\underline{v(g_{1}) \lhd
r(g_{2})\bigl) \lhd \bigl(v(g_{2}) \rhd r(g_{3})\bigl)\bigl]
\bigl(v(g_{2}) \lhd
r(g_{3}) \bigl)} v(g_{3})\\
&\stackrel{\lhd-{\rm~right~action~}}{=}& \bigl[v(g_{1}) \lhd
\underline{\bigl(r(g_{2})(v(g_{2})
\rhd r(g_{3}))\bigl)}\bigl]\bigl(v(g_{2}) \lhd r(g_{3})\bigl)v(g_{3})\\
&\stackrel{\equref{def1}} {=}& \bigl[v(g_{1}) \lhd r \circ
v^{-1}\bigl((v(g_{2}) \lhd r(g_{3}))v(g_{3})\bigl)\bigl]\\
&&\bigl(v(g_{2}) \lhd
r(g_{3})\bigl)v(g_{3})\\
&{=}& v\bigl[g_{1} * \underline{v^{-1}\bigl((v(g_{2}) \lhd r(g_{3})) v(g_{3})\bigl)} \bigl]\\
&\stackrel{\equref{def1}} {=}& v\bigl(g_{1} * (g_{2} *
g_{3})\bigl)
\end{eqnarray*}
i.e. the multiplication $*$ is associative as $v$ is a bijection.
Let $g \in G$ and define $g' := v^{-1}\bigl(v(g)^{-1} \lhd
r(g)^{-1}\bigl)$. Then:
\begin{eqnarray*}
v(g' * g) &=& \bigl( v(g') \lhd r(g) \bigl) v(g) = \Bigl( \bigl
(v(g)^{-1} \lhd r(g)^{-1} \bigl) \lhd r(g)\Bigl)
v(g)\\
&=& v(g)^{-1} v(g) = 1 = v (1)
\end{eqnarray*}
i.e. $g' * g = 1$ as $v$ is bijective. Thus every element $g \in
G$ has a left inverse, i.e. $(G,*)$ is a group.

(2) The proof can be done directly through a long computation but
we prefer the following approach: first we remark that the
defining relations \equref{def1}, \equref{def2}, \equref{def3} are
exactly the compatibility conditions \equref{p4}, \equref{p2}, and
respectively \equref{p1} from \prref{1} and the compatibility
condition \equref{crossed} is exactly \equref{p3} with the $*$
operations as defined by \equref{def1}. Moreover the map
$$
\psi
: H \times  (G,*) \rightarrow H\, {}_{\alpha}\!\! \bowtie_{\beta}
\, G, \quad \psi(h,g) = \bigl(\sigma(h)r(g), v(g)\bigl)
$$
is a bijection between the set $ H \times (G, *)$ and the group
$H\, {}_{\alpha}\!\! \bowtie_{\beta} \, G$. With this observation
in mind, in order to prove that $\bigl(H, (G,*), \alpha',
\beta'\bigl)$ is a matched pair it is enough to show that the
group structure obtained by transferring the group structure from
the bicrossed product $H\, {}_{\alpha}\!\! \bowtie_{\beta} \, G,$
to the set $ H \times (G, *)$ via the bijective map $\psi$ is
exactly the one of a bicrossed product on the set $ H \times (G,
*)$ associated to the actions $\alpha'$ and $\beta'$. In other
words, we have to prove that
$$
\bigl( h( g \rhd' h'), \, (g \lhd' h')g'\bigl) = \psi^{-1}
\Bigl(\psi(h,g) \cdot \psi(h',g')\Bigl)
$$
for all $h$, $h' \in H$, $g$, $g' \in G$ or equivalently, as
$\psi$ is bijective
\begin{equation}\eqlabel{def4}
\psi\bigl((h,g)\cdot (h',g')\bigl) = \psi(h,g) \cdot \psi(h',g')
\end{equation}
for all $h$, $h' \in H$, $g$, $g' \in G$. This reduces to proving
the following two conditions:
\begin{equation}\eqlabel{c11}
\sigma\bigl(h(g \rhd' h')\bigl) r\circ v^{-1}\bigl[\bigl(v(g \lhd'
h') \lhd r(g')\bigl)v(g')\bigl] = \sigma(h) r(g)\bigl(v(g) \rhd
(\sigma(h')r(g'))\bigl)
\end{equation}
and
\begin{equation}\eqlabel{c2}
\bigl(v(g \lhd' h') \lhd r(g')\bigl)v(g') = \bigl(v(g) \lhd
(\sigma(h')r(g'))\bigl)v(g')
\end{equation}
for any $h$, $h' \in H$, $g$, $g' \in G$. We have:
$$\bigl(v(g \lhd' h') \lhd r(g')\bigl) \stackrel{\equref{def2}} {=}
\bigl[\bigl(v(g) \lhd \sigma(h')\bigl) \lhd r(g')\bigl]v(g') =
\bigl(v(g) \lhd (\sigma(h')r(g'))\bigl)v(g')$$ Moreover:
\begin{eqnarray*}
\sigma(h) r(g)\bigl(v(g) \rhd
(\sigma(h')r(g'))\bigl)&\stackrel{\equref{2}} {=}& \sigma(h)
r(g)\bigl(v(g) \rhd \sigma(h')\bigl)\bigl((v (g) \lhd
\sigma(h')) \rhd r(g')\bigl)\\
&\stackrel{\equref{def2}} {=}& \sigma(h) r(g)\bigl(v(g) \rhd
\sigma(h')\bigl)\bigl(v (g \lhd'h') \rhd r(g')\bigl)\\
&{=}& \sigma(h) r(g)\bigl(v(g) \rhd \sigma(h')\bigl)r(g \lhd'
h')^{-1}r(g \lhd' h')\\
&&\bigl(v (g \lhd'h') \rhd r(g')\bigl)\\
&\stackrel{\equref{crossed}} {=}& \sigma(h)
r(g)\bigl(v(g) \rhd \sigma(h')\bigl) r(g \lhd' h')^{-1}r\bigl((g\lhd' h') * g'\bigl)\\
&\stackrel{\equref{def2},\equref{def1}} {=}& \sigma(h)
r(g)\bigl(v(g) \rhd \sigma(h')\bigl) r\circ v^{-1}\bigl(v(g) \lhd
\sigma(h)\bigl)^{-1}\\
&& r\circ v^{-1}\bigl[\bigl(v (g
\lhd'h')\lhd r(g')\bigl)v(g')\bigl]\\
&\stackrel{\equref{def3}} {=}& \sigma(h)\sigma(g \rhd' h') r\circ
v^{-1}\bigl[\bigl(v(g \lhd' h')\lhd r(g')\bigl)v(g')\bigl]\\
&{=}& \sigma\bigl(h(g \rhd' h')\bigl) r\circ v^{-1}\bigl[\bigl(v(g
\lhd' h') \lhd r(g')\bigl)v(g')\bigl]
\end{eqnarray*}
hence \equref{c11} holds. It follows that \equref{c2} and hence
\equref{def4} holds and we are done.

(3) Follows from (2) and the \prref{1}.
\end{proof}

\section{Schreier type theorems for bicrossed products}\selabel{3}

In this section we shall prove two Schreier type classification
theorems for bicrossed products. Let $H$ and $G$ be two fixed
groups.

Let $MP(H,G):= \{(\alpha,\beta) ~|~ (H, G, \alpha, \beta)
{\rm~is~}{\rm~a~}{\rm~matched~}{\rm~pair}\}$. We define
$B_{1}(H,G)$ to be the category having as objects the set
$MP(H,G)$ and the morphisms defined as follows: $\psi : (\alpha',
\beta') \to (\alpha, \beta)$ is a morphism in $B_{1}(H,G)$ if and
only if $\psi: H\, {}_{\alpha'}\!\! \bowtie_{\beta'} \, G
\rightarrow H\, {}_{\alpha}\!\! \bowtie_{\beta} \, G$ is a
morphism of groups such that $\psi \circ i_H = i_H$, where $i_H
(h) = (h, 1)$ is the canonical inclusion. Thus a morphism in the
category $B_{1}(H,G)$ is a morphism between two bicrossed products
of $H$ and $G$ that fix $H$. Considering $\sigma = Id_{H}$ and $G'
= G$ in \prref{1} we obtain the following:

\begin{corollary}\colabel{cosch}
Let $(H, G, \alpha, \beta)$ and $(H, G, \alpha', \beta')$ two
matched pairs. There exists a one to one correspondence between
the set of all morphisms $\psi: H\, {}_{\alpha'}\!\!
\bowtie_{\beta'} \, G \rightarrow H\, {}_{\alpha}\!\!
\bowtie_{\beta} \, G$ in the category $B_{1}(H,G)$ and the set of
all pairs $(r, v)$, where $r: G \rightarrow H$, $v: G \rightarrow
G$ are two maps such that:
\begin{eqnarray}
(g \rhd' h)r(g \lhd' h) &=& r(g)(v(g) \rhd
h)\eqlabel{s1} \\
v(g \lhd' h) &=& v(g) \lhd h \eqlabel{s2} \\
r(g_{1}g_{2}) &=& r(g_{1})\bigl(v(g_{1}) \rhd
r(g_{2}) \bigl) \eqlabel{s3} \\
v(g_{1}g_{2}) &=& \bigl(v(g_{1}) \lhd r(g_{2})\bigl)v(g_{2})
\eqlabel{s4}
\end{eqnarray}
for all $h \in H$, $g$, $g_{1}$, $g_{2} \in G$. Through the above
bijection $\psi$ is given by
\begin{equation}\eqlabel{d}
\psi(h,g) = \bigl(hr(g),v(g)\bigl)
\end{equation}
and $\psi : H\, {}_{\alpha'}\!\! \bowtie_{\beta'} \, G \rightarrow
H\, {}_{\alpha}\!\! \bowtie_{\beta} \, G$ is an isomorphism of
groups that fixes $H$ if and only if $v: G\to G$ is a bijective
map.
\end{corollary}

We are led to the following:

\begin{definition}
Let $H$ and $G$ be two groups. Two pairs $(\alpha, \beta)$,
$(\alpha', \beta') \in MP(H,G)$ are called \textit{1-equivalent}
and we denote this by $(\alpha, \beta) \approx_{1} (\alpha',
\beta')$ if and only if there exists a pair $(r, v)$, where $r: G
\rightarrow H$, $v: G \rightarrow G$ are two maps such that $v$ is
bijective and the relations \equref{s1} - \equref{s4} hold.
\end{definition}

Thus, using \coref{cosch}, we obtain that $(\alpha, \beta)
\approx_{1} (\alpha', \beta')$ if and only if there exists an
isomorphism $(\alpha, \beta) \cong (\alpha', \beta')$ in the
category $B_{1}(H,G)$. In particular, $\approx_{1}$ is an
equivalence relation on the set $MP(H,G)$ and we have proved:

\begin{theorem}\textbf{(First Schreier type theorem for bicrossed
products)}\thlabel{sch1} Let $H$ and $G$ be two groups. There
exists a bijection between the set of objects of the skeleton of
the category $B_{1}(H,G)$ and the pointed quotient set
$MP(H,G)/\approx_{1}$. We shall use the following notation:
$K^{2}(H,G):= M(H,G)/\approx_{1}$.
\end{theorem}

A general problem arises in order to have a classification type
theorem for bicrossed products of two given groups $H$ and $G$:

\textit{Compute $K^{2}(H,G)$ for two given groups $H$ and $G$.}

In the last section we shall compute explicitly the set $K^2 (C_3,
C_6)$. $K^{2}(H,G)$ is a pointed set by the equivalence class of
the pair $(\alpha_{0}$, $\beta_{0})$, where $\alpha_{0}$,
$\beta_{0}$ are the trivial actions. For $\alpha := \alpha_{0}$
and $\beta := \beta_{0}$ we obtain from relations \equref{s1} -
\equref{s4} that $\lhd'$ is the trivial action, $r: G\to H$ and
$v: G\to G$ are morphisms of groups and $g \rhd' h =
r(g)hr(g)^{-1}$, for all $g\in G$ and $h\in H$. For every morphism
of groups $r: G\to H$ we shall denote by $\rhd_{r}$ the action $g
\rhd_{r} h := r(g)hr(g)^{-1}$, for all $g\in G$ and $h\in H$.
Hence $\widehat{(\alpha_{0}, \beta_{0})} = \{(\rhd_{r}, \beta_{0})
~|~ r: G\to H ~~ {\rm ~is~ a~ morphism~ of~ groups}\}$. We restate
this observation as follows: let $H$ and $G$ be two groups. Then
there exists $(H, G, \alpha', \beta')$ a matched pair such that
$H\, {}_{\alpha'}\!\! \bowtie_{\beta'} \, G \cong H\times G$
(isomorphism of groups that fixes $H$) if and only if the action
$\beta'$ is trivial and there exists a morphism of groups $r: G
\rightarrow H$ such that the action $\alpha'$ is given by $g \rhd'
h = r(g) h r(g)^{-1}$ for all $g \in G$, $h \in H$.

More generally, as a first application of \thref{sch1}, we shall
prove the following necessary and sufficient condition for a
bicrossed product to be isomorphic to a left version of a
semidirect product in the category $B_{1}(H,G)$:

\begin{corollary}\colabel{lake1}
Let $H$, $G$ be two groups and $\alpha : G\times H \to H$ be an
action as automorphisms of $G$ on $H$. The following statements
are equivalent:
\begin{enumerate}
\item There exists $(H, G, \alpha', \beta')$ a matched pair of
groups such that $H\, {}_{\alpha'}\!\! \bowtie_{\beta'} \, G \cong
H {}_{\alpha}\ltimes  G$ an isomorphism of groups that fixes $H$.

\item The action $\beta'$ is trivial and there exists a pair $(r,
v)$, where $v\in \Aut (G)$ is an automorphism of $G$, $r: G \to H$
is a map such that
\begin{equation}\eqlabel{s3'}
r(g_{1}g_{2}) = r(g_{1})\bigl(v(g_{1}) \rhd r(g_{2}) \bigl)
\end{equation}
for all $g_1$, $g_2 \in G$ and the action $\alpha '$ is
implemented as follows:
\begin{equation}\eqlabel{s1'}
g \rhd' h = r(g)(v(g) \rhd h) r(g)^{-1}
\end{equation}
for all $g\in G$ and $h\in H$.
\end{enumerate}
The isomorphism $\psi: H\, {}_{\alpha'}\!\! \bowtie_{\beta'} \, G
\to H {}_{\alpha}\ltimes G$ in $B_{1}(H,G)$ is given by $\psi (h,
g) = (h r(g), v(g))$, for all $h\in H$, $g\in G$.
\end{corollary}

\begin{proof} We apply \coref{cosch} in the case that $\beta$ is
the trivial action. It this context, using the fact that $v$ is
bijective, it follows from \equref{s2} that the action $\beta'$ is
trivial and \equref{s4} reduces to the fact that $v$ is a
morphism, hence an automorphism of $G$. Finally, \equref{s1} and
\equref{s3} are exactly \equref{s1'} and \equref{s3'}.
\end{proof}

We shall give now a necessary and sufficient condition for a
bicrossed product to be isomorphic to a right version of a
semidirect product in the category $B_{1}(H,G)$:

\begin{corollary}\colabel{lake2}
Let $H$, $G$ be two groups and $\beta : G\times H \to G$ be an
action as automorphisms of $H$ on $G$. The following statements
are equivalent:
\begin{enumerate}
\item There exists $(H, G, \alpha', \beta')$ a matched pair of
groups such that $H\, {}_{\alpha'}\!\! \bowtie_{\beta'} \, G \cong
H \rtimes_{\beta} G$ is an isomorphism of groups that fixes $H$.

\item There exists a pair $(r, v)$, where $r: G \to H$ is a
morphism of groups, $v: G\to G$ is a bijective map such that

\begin{equation}\eqlabel{s4''}
v(g_{1}g_{2}) = \bigl(v(g_{1}) \lhd r(g_{2})\bigl)v(g_{2})
\end{equation}
for all $g_1$, $g_2 \in G$ and the actions $\alpha '$ and $\beta'$
are implemented as follows:
\begin{eqnarray}
g \rhd' h &=& r(g) \, h \, \Bigl ( r\circ v^{-1} \bigl (v(g) \lhd
h\bigl)
\Bigl)^{-1} \eqlabel{s1''}\\
g \lhd' h &=& v^{-1} \bigl( v(g) \lhd h \bigl) \eqlabel{s2''}
\end{eqnarray}
for all $g\in G$ and $h\in H$.
\end{enumerate}
The isomorphism $\psi: H\, {}_{\alpha'}\!\! \bowtie_{\beta'} \, G
\to H \rtimes_{\beta} G$ in $B_{1}(H,G)$ is given by $\psi (h, g)
= (h r(g), v(g))$, for all $h\in H$, $g\in G$.
\end{corollary}

\begin{proof} We apply \coref{cosch} in the case that $\alpha$ is
the trivial action. It this context, it follows from \equref{s3}
that $r$ is a morphism of groups, while \equref{s1} and
\equref{s2} are exactly \equref{s1''} and \equref{s2''}.
\end{proof}

In what follows we will prove the second Schreier type theorem for
bicrossed products: it is the analogue of the theorem regarding
group extensions. Let $(H, G, \alpha, \beta)$ be a matched pair.
Then the natural projections $\pi_{G}:H\, {}_{\alpha}\!\!
\bowtie_{\beta} \, G \rightarrow G$, $\pi_{H}:H\, {}_{\alpha}\!\!
\bowtie_{\beta} \, G \rightarrow H$ are not morphisms of groups.

We will fix two groups $H$, $G$ and $\beta: G \times H \rightarrow
G$ a right action of the group $H$ on the set $G$. We denote by
$\widetilde{\beta}: H \rightarrow End(G)$ the corresponding
morphism of groups. Define
$${\rm Ker}(\widetilde{\beta}):= \{h \in H ~|~ g \lhd h = g, \forall g \in G\}$$
We denote by $MP_{\beta}(H,G):= \{\alpha ~|~ (H, G, \alpha, \beta)
{\rm~is~}{\rm~a~}{\rm~matched~}{\rm~pair~}\}$. Let
$B_{2}^{\beta}(H,G)$ be the category having $MP_{\beta}(H,G)$ as
the set of objects and the morphisms defined as follows: $\psi:
\alpha' \rightarrow \alpha$ is a morphism in $B_{2}^{\beta}(H,G)$
if and only if $\psi: H\, {}_{\alpha'}\!\! \bowtie_{\beta} \, G
\rightarrow H\, {}_{\alpha}\!\! \bowtie_{\beta} \, G$ is a
morphism of groups such that
\begin{equation}\eqlabel{diag2}
\psi \circ i_H = i_H ~~~{\rm and}~~~ \pi_G \circ \psi = \pi_G
\end{equation}

\begin{proposition}\prlabel{2}
Let $(H, G, \alpha', \beta)$, $(H, G, \alpha, \beta)$ be two
matched pairs. There exists a one to one correspondence between
the set of all morphisms $\psi: \alpha' \rightarrow \alpha$ in the
category $B_{2}^{\beta}(H,G)$ and the set of all maps $r:G
\rightarrow {\rm Ker} (\widetilde{\beta})$ such that :
\begin{eqnarray}
(g \rhd' h)r(g \lhd h) &=& r(g)(g \rhd
h)\eqlabel{p1'} \\
r(g_{1} g_{2}) &=& r(g_{1})\bigl(g_{1} \rhd r(g_{2})\bigl)
\eqlabel{p3'}
\end{eqnarray}
for all $g, g_{1}, g_{2} \in G$, $h \in H$. Through the above
bijection the morphism $\psi$ is given by
\begin{equation}\eqlabel{p5'}
\psi(h,g) = \bigl(hr(g), g\bigl)
\end{equation}
for all $h \in H$, $g \in G$ and $\psi$ is an isomorphism of
groups i.e. $B_{2}^{\beta}(H,G)$ is a groupoid \footnote {We
recall that a groupoid is a category such that the class of
objects is a set and any morphism is an isomorphism.}.
\end{proposition}

\begin{proof}
For any morphism of groups $\psi: H\, {}_{\alpha'}\!\!
\bowtie_{\beta} \, G \rightarrow H\, {}_{\alpha}\!\!
\bowtie_{\beta} \, G$ such that \equref{diag2} hold there exists a
unique map $r: G\to H$ such that $\psi (h, g) = (h r(g), g)$, for
all $h\in H$ and $g\in G$. Now we are in a position to use
\prref{1} for $G' = G$, $\beta' = \beta$ and $v = Id_G$. We obtain
\equref{p1'} and \equref{p3'} by considering $v = Id_{G}$ in
\equref{p1}, respectively \equref{p3}. On the other hand
\equref{p2} is trivially fulfilled and \equref{p4} becomes $g_{1}
= g_{1} \lhd r(g_{2})$ for all $g_{1}$, $g_{2} \in G$, i.e. ${\rm
Im}(r) \subseteq {\rm Ker}(\widetilde{\beta})$.
\end{proof}

\begin{remark}
If $\beta$ is a faithful action (i.e. ${\rm
Ker}(\widetilde{\beta}) = 1$) then $B_{2}^{\beta}(H,G)$ is a
discret groupoid i.e. there exists a morphism $\psi: \alpha'
\rightarrow \alpha$ if and only if $\alpha = \alpha'$. Indeed, in
this case $r(g) = 1$ for all $g \in G$, \equref{p3'} is trivially
fulfilled and \equref{p1'} reduces to $g \rhd h = g \rhd' h$ i.e.
$\alpha = \alpha'$. Hence, in this case the skeleton of the
category $B_{2}^{\beta}(H,G)$ is the set $MP_{\beta}(H,G)$.
\end{remark}

\begin{definition}
Let $H$, $G$ be two groups and $\beta : G \times H \rightarrow G$
be a right action. Two elements $\alpha'$ and $\alpha$ of
$MP_{\beta}(H,G)$ are called $\approx_{2}$-equivalent and we
denote this by $\alpha' \approx_{2} \alpha$ if there exists a map
$r:G \rightarrow {\rm Ker}(\widetilde{\beta})$ such that the
relations \equref{p1'} and \equref{p3'} hold.
\end{definition}

From \prref{2} we obtain that $\alpha' \approx_{2} \alpha$ if and
only if there exists an isomorphism $\alpha' \cong \alpha$ in
$B_{2}^{\beta}(H,G)$. Hence $\approx_{2}$ is an equivalence
relation on $MP_{\beta}(H,G)$ and we obtained the following:

\begin{theorem}\textbf{(The second Schreier type theorem for
bicrossed products)}\thlabel{sch22} Let $H$, $G$ be two groups and
$\beta: G \times H \rightarrow G$ be a right action. There exists
a bijection between the set of objects of the skeleton of the
category $B_{2}^{\beta}(H,G)$ and the quotient set
$MP_{\beta}(H,G)/\approx_{2}$.
\end{theorem}

It is possible that the set $MP_{\beta}(H,G)$ (and hence
$MP_{\beta}(H,G)/\approx_{2}$) is the empty set. However, if
$\beta : G\times H \to G$ is an action as automorphisms then
$MP_{\beta}(H,G)$ is nonempty as it contains the trivial action
$\alpha_0$. In this case the quotient set
$MP_{\beta}(H,G)/\approx_{2}$ is a pointed set by the equivalence
class of the trivial action $\alpha_0$. It follows from \prref{2}
that
$$\widehat{\alpha_{0}} = \{\alpha' | \, \alpha' (g, h) = r(g) h
r(g \lhd h)^{-1}, {\rm for ~ some ~} r: G\to {\rm
Ker}(\widetilde{\beta}) ~~ {\rm ~ a~ morphism~ of~ groups}\}$$

We record this observation in the following:

\begin{corollary}
Let $H$, $G$ be two groups, $\beta: G \times H \rightarrow G$ an
action as automorphisms and $H\rtimes_{\beta} G$ the right version
of the semidirect product. The following statements are
equivalent:
\begin{enumerate}
\item There exists a matched pair $(H, G, \alpha, \beta)$ such
that the bicrossed products $H\, {}_{\alpha}\!\! \bowtie_{\beta}
\, G$ and $H\rtimes_{\beta} G$ are isomorphic in the category
$B_{2}^{\beta}(H,G)$;

\item There exists a morphism of groups $r: G \rightarrow {\rm
Ker}(\widetilde{\beta})$ such that the action $\alpha$ is given by
$g \rhd h = r(g) h r(g \lhd h)^{-1}$ for all $g \in G$, $h \in H$.
\end{enumerate}
\end{corollary}

\section{Examples}\selabel{4}

In this section we describe all matched pairs between $C_{n}$ and
$C_{m}$, for $n \in \{2,3\}$ and $m \in \NN^{*}$ arbitrary. First,
let us introduce some notation. We denote by $a$ a generator of
the cyclic group $C_n$ and $b$ a generator of $C_m$. The set of
group  morphisms from the group $C_n$ to the group of
automorphisms $\Aut (C_m)$ will be denoted by $\varsigma (n, m)$.
Such a morphism $\vartheta : C_n \rightarrow \Aut (C_m)$ is
uniquely determined by a positive integer $t\in [m-1] := \{1, 2,
\cdots, m-1\} $ such that $m|t^n -1$ and
\begin{equation}\eqlabel{2.4.399}
\vartheta : C_n \rightarrow \Aut (C_m), \qquad \vartheta (a) (b) =
b^t
\end{equation}
Therefore, one can equivalently think of $\varsigma (n, m)$ as the
subgroup of $U(\ZZ_m)$ consisting of all solutions in $\ZZ_{m}$ of
the equation $x^n = 1$.

Using the fact that if $m= 2^{a_0} p_1^{a_1}\cdots p_k^{a_k}$ with
$p_1$, $\cdots$, $p_k$  odd primes, then  $$\Aut (C_m)\cong U(
\ZZ_m ) \cong U( \ZZ_{2^{a_0}}) \times U(\ZZ_{p_1^{a_1}}) \times
\cdots \times U(\ZZ_{p_k^{a_k}})$$ it is a routine computation to
check that $$ |\varsigma (n, m)| = \left \{\begin{array}{rcl}

\prod_{i=1}^k (n, p_i^{a_i} - p_i^{a_i -1}), \, & \mbox { if }&
4\!\not|
m\\
(n,2) (n, 2^{a_0 -2}) \prod_{i=1}^k (n, p_i^{a_i} - p_i^{a_i -1}),
\, & \mbox { if }& 4|m
\end{array} \right.
$$
In particular
$$
|\varsigma (2, m) | = \left \{\begin{array}{rcl}
2^k,  \, & \mbox {\rm if }& \, a_0\leq 1\\
2^{k+1},  \, & \mbox {\rm if }& \, a_0=2\\
2^{k+2},  \, & \mbox {\rm if }& \, a_0\geq 3
\end{array} \right.
$$
and
$$
| \varsigma (p, m) | = \left \{\begin{array}{rcl}
\prod_{i=1}^k (p, p_i -1),  \, & \mbox {\rm if }& \, p^2\!\not| m\\
p\prod_{i=1}^k (p, p_i -1),  \, & \mbox {\rm if }& \, p^2|m
\end{array} \right.
$$
for an odd prime $p$.

Let $m$ be a positive integer. Then $(C_2, C_m, \alpha, \beta)$ is
a matched pair if and only if the action $\alpha$ is trivial and
there exists a positive integer $t\in [m-1]$ such that $m|t^2 -1$
and $\beta = \beta_t : C_m \times C_2 \rightarrow C_m$ is given by
\begin{equation}\eqlabel{2.4.6}
\beta (b^i, a) = b^{it}, \quad \beta (b^i, 1) = b^i
\end{equation}
for any $i= 0, \cdots, m-1$. In particular, there are $|\varsigma
(2, m)|$ matched pairs $(C_2, C_m, \alpha, \beta)$.

Indeed, as $\alpha$ is an action we get $b\triangleright a \neq 1
= b \triangleright 1$. Thus $b \triangleright a = a$ which implies
that $\alpha$ is trivial. Thus $(C_2, C_m, \alpha, \beta)$ is a
matched pair if and only if $\beta ' : C_2 \rightarrow \Aut
(C_m)$, $\beta' (x) (y) := \beta (y, x)$ is a morphism of groups,
so by letting $n=2$ in \equref{2.4.399} we obtain that $(C_2, C_m,
\alpha, \beta)$ is a matched pair if and only if there exists
$t\in [m-1]$ such that $m|t^2 -1$ and $\beta (b, a) = b^t$. The
formula \equref{2.4.6} follows as $\beta$ is an action.

In order to describe all matched pairs $(C_3, C_m, \alpha, \beta)$
we need the following observation.

\begin{remark}\relabel{knitauto}
Let $(H, G, \alpha, \beta)$ be a matched pair such that $\alpha$
is an action of $G$  on $H$ as group automorphisms. Then the
compatibility condition \equref{2} from the definition of a
matched pair is equivalent to $(g\triangleleft {h_1})
\triangleright h_2 = g \triangleright h_2$ that can be written as
\begin{equation}\eqlabel{KS4'}
g^{-1}(g\triangleleft {h_1}) \in {\rm Stab}_G (h_2)
\end{equation}
for any $g\in G$, $h_1$, $h_2\in H$. Thus if $\alpha$ is an action
as automorphisms then $(H, G, \alpha, \beta)$ is a matched pair if
and only if \equref{3} and \equref{KS4'} hold. The condition
\equref{KS4'} gives  important information regarding  $\beta$: the
elements $g^{-1}\beta (g, h)$ act trivially on $H$ for any $g\in
G$ and $h\in H$.
\end{remark}

Now we can describe all matched pairs $(C_3, C_m, \alpha, \beta)$.

\begin{proposition}\prlabel{2.4.45}
Let $m$ be a positive integer, $\alpha : C_m \times C_3
\rightarrow C_3$, $\beta : C_m \times C_3 \rightarrow C_m$ two
maps and $t\in [m-1]$ such that $m|t^3 -1$. Then:
\begin{enumerate}

\item[(i)] Let $\alpha$ be the trivial action and $\beta = \beta_t
: C_m \times C_3 \rightarrow C_m$ given by
\begin{equation}\eqlabel{2.4.90}
\beta (b^i, a) = b^{it}, \quad  \beta (b^i, a^2) = b^{it^2}, \quad
\beta (b^i, 1) = b^i
\end{equation}
for any $i= 0, \cdots, m-1$. Then $(C_3, C_m, \alpha, \beta_t)$ is
a matched pair. There are no other matched pairs $(C_3, C_m,
\alpha, \beta)$ if $m$ is odd.

\item[(ii)] Assume that $m$ is even. Let $\beta$ be the trivial
action and $\alpha : C_m \times C_3 \rightarrow C_3$ given by
$\alpha (b^j, 1) = 1$ and
\begin{equation}\eqlabel{2.4.740}
\alpha (b^j, a) = \left \{\begin{array}{rcl}
a,  \, & \mbox {\rm if $j$ is even }\\
a^2,  \, & \mbox {\rm if $j$ is odd }
\end{array} \right.
\end{equation}
\begin{equation}\eqlabel{2.4.750}
\alpha (b^j, a^2) = \left \{\begin{array}{rcl}
a^2,  \, & \mbox {\rm if $j$ is even }\\
a,  \, & \mbox {\rm if $j$ is odd }
\end{array} \right.
\end{equation}
for all $j = 1, \cdots, m-1$. Then $(C_3, C_m, \alpha, \beta)$ is
a matched pair.

\item[(iii)] Assume that $m = 6 u$ for some positive integer $u$
and that $\alpha$ is described by \equref{2.4.740} and
\equref{2.4.750}. Then there exist two matched pairs $(C_3, C_m,
\alpha, \beta)$, $(C_3, C_m, \alpha, \beta')$, where $\beta$ and
$\beta'$ are given by
\begin{equation}\eqlabel{24.11}
\beta (b^{2k+1}, a) = b^{2u +2k +1}, \qquad \beta (b^{2k+1}, a^2)
= b^{4u+ 2k +1}
\end{equation}
\begin{equation}\eqlabel{24.11'}
\beta (b^{2k}, a) = \beta (b^{2k}, a^2) = b^{2k}
\end{equation}
and
\begin{equation}\eqlabel{24.12}
\beta' (b^{2k+1}, a) = b^{4u +2k +1}, \qquad \beta' (b^{2k+1},
a^2) = b^{2u+ 2k + 1}
\end{equation}
\begin{equation}\eqlabel{24.12'}
\beta'(b^{2k}, a) = \beta'(b^{2k}, a^2) = b^{2k}
\end{equation}
for all nonnegative integers $k$. In this case there are $2+
|\varsigma (3, m)|$ matched pairs between $C_3$ and $C_m$.
\item[(iv)] There are no other matched pairs on $(C_3, C_m,
\alpha, \beta)$ other than the ones described above.
\end{enumerate}
\end{proposition}

\begin{proof} We assume first that $\alpha$ is the trivial action.
It follows from \reref{2.4.90} that $(C_3, C_m, \alpha, \beta)$ is
a matched pair if and only if $\beta ' : C_3 \rightarrow \Aut
(C_m)$, $\beta' (x) (y) := \beta (y, x)$ is a morphism of groups;
setting $n=3$ in \equref{2.4.399} we get that $(C_3, C_m, \alpha,
\beta)$ is a matched pair if and only if there exists $t\in [m-1]$
such that $m|t^3 -1$ and $\beta (b, a) = b^t$. The formula
\equref{2.4.90} follows as $\beta$ is an action.

Assume now that $m$ is odd. It follows from \equref{4} that $b^i
\triangleright 1 = 1$ for all $i= 0, \cdots, m-1$. If $b
\triangleright a = a$ we obtain that  $\alpha$ is trivial. Assume
that $b \triangleright a = a^2$. Then $b \triangleright a^2 = a$
and $\alpha$ is given by \equref{2.4.740}, \equref{2.4.750}. If
$m$ is odd we obtain: $a = 1\triangleright a = b^m \triangleright
a = a^2$, contradiction. Hence, for an odd $m$ the action $\alpha$
must be trivial and (i) is proved.

Assume now that $m$ is even and $\beta$ is the trivial action.
Then $(C_3, C_m, \alpha, \beta)$ is a matched pair if and only if
$\alpha$ is an action of $C_m$ on $C_3$ as group automorphisms.
The map $\alpha$ given by \equref{2.4.740}, \equref{2.4.750} is
such an action corresponding to
$$C_m \rightarrow \Aut (C_3), \quad b \mapsto (a \mapsto a^2)$$
Therefore (ii) is proved.

We shall prove now (iii) and (iv). Let $(C_3, C_m, \alpha, \beta)$
be a matched pair. We have proved that $\alpha : C_m \times C_3
\rightarrow C_3$ is either the trivial action or it is given by
\equref{2.4.740}, \equref{2.4.750} if $m$ is even.

We assume now that $m$ is even and that $\alpha$ is given by
\equref{2.4.740}, \equref{2.4.750}. Then $\alpha$ is an action of
$C_m$ on $C_3$ as automorphisms and
$${\rm Stab}_{C_m} (a) = {\rm
Stab}_{C_m} (a^2) = <b^2>$$ Using \reref{knitauto} we get that
$b^{-1}\beta (b, a) \in \, <b^2>$ and $b^{-1}\beta (b, a^2) \in \,
<b^2>$ and \equref{2} holds automatically. Let $l$, $t \in \{0, 1,
\cdots, m/2 -1 \}$ such that
\begin{equation}\eqlabel{4.444}
\beta (b, a) = b^{2l+1}, \qquad \beta (b, a^2) = b^{2t+1}
\end{equation}
We shall extend $\beta$ for each element of $C_m \times C_3$ using
\equref{3} as defining relations and the fact that $\beta$ is an
action. First we define $\beta$ for each pair $(b^i, a)$ such that
\equref{3} holds. We have
$$
\beta (b^2, a ) = (bb)\triangleleft a \stackrel{\equref{3} } =
(b\triangleleft {(b\triangleright a)}) (b\triangleleft a)
\stackrel{\equref{2.4.740} } = (b\triangleleft {a^2})(b
\triangleleft a) = b ^{2(t+l +1)}
$$
and
$$
\beta (b^3, a ) = (bb^2)\triangleleft a \stackrel{\equref{3} } =
(b\triangleleft {(b^2\triangleright a)}) (b^2\triangleleft a)
\stackrel{\equref{2.4.740} } = (b\triangleleft {a})
(b^2\triangleleft a) = b ^{4l + 2t +3}
$$
Using the induction we can prove
\begin{equation}\eqlabel{4.445}
\beta (b^{2k}, a) =  b^{2k(l+t+1)} , \qquad \beta (b^{2k+1}, a) =
b^{(2k+2)l+2kt+2k+1}
\end{equation}
for any $k = 0, 1, \cdots$. We note that
$$
b^{2t+1} = \beta (b, a^2) = b\triangleleft {a^2} = (b\triangleleft
a)\triangleleft a = \beta (b^{2l+1}, a)\stackrel{\equref{4.445} }
= b^{(2l+2)l + 2lt + 2l + 1}
$$
As the order of $b$ is $m$ we get a first compatibility condition
for $l$ and $t$:
\begin{equation}\eqlabel{c1}
m | 2(l^2 + 2l + lt -t)
\end{equation}
Now we define $\beta$ for each pair $(b^i, a^2)$ using \equref{3}
repeatedly. We have:
$$
\beta (b^2, a^2) = (bb)\triangleleft {a^2} \stackrel{\equref{3} }=
(b\triangleleft {(b\triangleright a^2)}) (b\triangleleft {a^2})
\stackrel{ \equref{2.4.750} } =( b\triangleleft a)( b\triangleleft
{a^2}) = b^{2(l + t + 1 )}
$$
$$
\beta (b^3, a^2) = (bb^2)\triangleleft {a^2} \stackrel{\equref{3}
}= (b\triangleleft {(b^2\triangleright a^2)}) (b^2 \triangleleft
{a^2}) \stackrel{ \equref{2.4.750} } = (b\triangleleft {a^2})
(b^2\triangleleft {a^2}) = b^{2l + 4t + 3}
$$
Using the induction we can easily prove that
\begin{equation}\eqlabel{4.445'}
\beta (b^{2k}, a^2) =  b^{2k(l+t+1)} , \qquad \beta (b^{2k+1},
a^2) = b^{ 2kl + (2k+2)t + 2k + 1}
\end{equation}
for any $k = 0, 1, \cdots$. Moreover, keeping in mind
\equref{4.445} we find that
$$
\beta (b^{2k}, a) = \beta (b^{2k}, a^2) = b^{2k(l+t+1)}
$$
On the other hand $\beta$ is a right action and $a^3 = 1$. Hence:
$$
b^{2k} = \beta (b^{2k}, 1) = (b^{2k}\triangleleft
{a^2})\triangleleft a = (b^{2k}\triangleleft {a})\triangleleft a =
b^{2k}\triangleleft {a^2} =  b^{2k(l+t+1)}
$$
As the order of $b$ is $m$ we obtain a second compatibility
condition between $l$ and $t$: $m | 2k(l+t)$ for any $k = 0, 1,
\cdots$ which is equivalent to:
\begin{equation}\eqlabel{c5}
m | 2(l+t)
\end{equation}
From this condition and \equref{c1} we obtain
\begin{equation}\eqlabel{c1'}
m | 2(2l - t)
\end{equation}
Let now $m = 2r$. We have to find $l$, $t \in \{1, 2, \cdots, r-1
\}$ such that
$$
m | 2(l+t) \quad {\rm and} \quad m | 2(2l -t)
$$
Equivalently,  we have to solve in  $\ZZ_r$ the system of
equations
\begin{equation}\eqlabel{4.600}
\hat{l} + \hat{t} = \hat{0}, \qquad \hat{2}\hat{l} - \hat{t} =
\hat{0}
\end{equation}
The equation $\hat{3} \hat{l} = \hat{0}$ has $(3, r)$ solutions in
$\ZZ_r$. If $3$ does not divide  $m$ then the unique solution of
the system is $\hat{l} = \hat{t} = \hat{0}$ and  therefore $\beta$
is the trivial action. If $3$ divides $r$ let $u$ be such that $r
= 3u$. Then the system \equref{4.600} has three solutions
$$\hat{l_1} =
\hat{t}_1 = \hat{0}; \qquad \hat{l}_2 = \hat{u}, \, \hat{t}_2 =
\hat{2}\hat{u}, \qquad \hat{l_3} = \hat{2}\hat{u}, \, \hat{t_3} =
\hat{4}\hat{u}
$$
The first solution gives that the action $\beta$ is trivial and
the last two solutions give exactly the two actions $\beta$
described in \equref{24.11} and \equref{24.12}.
\end{proof}

We showed that the smallest example of a proper matched pair (i.e.
one in which both actions are nontrivial) between two finite
cyclic groups is the one between the groups $C_3$ and $C_6$.
According to \prref{2.4.45} there exist exactly four matched pairs
$(C_3, C_6, \alpha, \beta)$ namely:
\begin{enumerate}
\item[(i)] $\alpha_{0}$ and $\beta_{0}$ are the trivial actions;
\item[(ii)] $\beta_{0}$ is the trivial action and $\alpha_{1}$ is
defined by:
\begin{eqnarray*}
b^j\rhd_{1} a = \left \{\begin{array}{rcl}
a,  \, & \mbox {\rm if $j$ is even }\\
a^2,  \, & \mbox {\rm if $j$ is odd }
\end{array} \right.
, \quad b^j \rhd_{1} a^2 = \left \{\begin{array}{rcl}
a^2,  \, & \mbox {\rm if $j$ is even }\\
a,  \, & \mbox {\rm if $j$ is odd }
\end{array} \right.
\end{eqnarray*}
for all $j = 1, \cdots, 5$.
\item[(iii)] $\alpha_{2}$ and
$\beta_{2}$ are defined by :
\begin{eqnarray*}
b^j\rhd_{2} a = \left \{\begin{array}{rcl}
a,  \, & \mbox {\rm if $j$ is even }\\
a^2,  \, & \mbox {\rm if $j$ is odd }
\end{array} \right.
, \quad b^j \rhd_{2} a^2 = \left \{\begin{array}{rcl}
a^2,  \, & \mbox {\rm if $j$ is even }\\
a,  \, & \mbox {\rm if $j$ is odd }
\end{array} \right.
\end{eqnarray*}
and
\begin{eqnarray*}
b^j\lhd_{2} a = \left \{\begin{array}{rcl}
b^{j},  \, & \mbox {\rm if $j$ is even }\\
b^{j+2},  \, & \mbox {\rm if $j$ is odd }
\end{array} \right.
, \quad b^j \lhd_{2} a^2 = \left \{\begin{array}{rcl}
b^{j},  \, & \mbox {\rm if $j$ is even }\\
b^{j+4},  \, & \mbox {\rm if $j$ is odd }
\end{array} \right.
\end{eqnarray*}
for all $j = 1, \cdots, 5$.
\item[(iv)] $\alpha_{2}$ and
$\beta_{3}$ are defined by:
\begin{eqnarray*}
b^j\lhd_{3} a = \left \{\begin{array}{rcl}
b^{j},  \, & \mbox {\rm if $j$ is even }\\
b^{j+4},  \, & \mbox {\rm if $j$ is odd }
\end{array} \right.
, \quad b^j \lhd_{3} a^2 = \left \{\begin{array}{rcl}
b^{j},  \, & \mbox {\rm if $j$ is even }\\
b^{j+2},  \, & \mbox {\rm if $j$ is odd }
\end{array} \right.
\end{eqnarray*}
for all $j = 1, \cdots, 5$.
\end{enumerate}

We shall now classify all bicrossed products $C_3 \bowtie C_6$
that fixes the group $C_3$, i.e. we shall determine the pointed
set $K^{2}(C_{3},C_{6})$ from \thref{sch1}.

\begin{corollary}
$K^{2}(C_{3},C_{6})$ is a pointed set with three elements. In
particular, any bicrossed product $C_3 \bowtie C_6$ that fixes the
group $C_3$ is isomorphic to one of the following three groups:
$$
C_3 \times C_6, \quad < a, b \, | \, a^3 = 1, \, b^6 = 1, \, ba =
a^2 b > , \quad < a, b \,|\, a^3 = 1, \, b^6 = 1, \, ba = a^2b^3 >
$$
\end{corollary}

\begin{proof}
Let $(\alpha', \beta')$ be a matched pair such that $(\alpha_{0},
\beta_{0}) \approx_{1} (\alpha', \beta')$. The relations
\equref{s1} - \equref{s2} collapse into:
\begin{eqnarray}
(g \rhd' h) r(g \lhd ' h) &=& r(g)h
\eqlabel{ex1}\\
v(g \lhd ' h) &=& v(g) \eqlabel{ex2}
\end{eqnarray}
Since $v$ is a bijective map, it follows from \equref{ex2} that
$\beta'$ is the trivial action. Furthermore, from \equref{ex1} we
obtain that $\alpha'$ is also the trivial action, that is, the
equivalence class of $(\alpha_{0}, \beta_{0})$ is trivial. By
similar arguments it follows that the equivalence class of
$(\alpha_{1}, \beta_{0})$ is also trivial.

Consider now $r: C_{6} \rightarrow C_{3}$ be the trivial morphism
of groups and $v: C_{6} \rightarrow C_{6}$ the automorphism given
by $v(b) := b^{5}$. By a straightforward computation it follows
that:
\begin{eqnarray*}
g \rhd_{2} h = v(g)\rhd_{2} h, \quad v(g \lhd_{3} h) = v(g)
\lhd_{2} h
\end{eqnarray*}
hence $(\alpha_{2}, \beta_{2}) \approx_{1} (\alpha_{2},
\beta_{3})$. Thus, $K^{2}(C_{3},C_{6})$ is a set with three
elements.
\end{proof}

\begin{remark}
It is easy to see that $B_{2}^{\beta}(C_{3},C_{6})$ is a singleton
or a set with two elements for any right action $\beta$.

Indeed, it is obvious that $B_{2}^{\beta_{2}}(C_{3},C_{6}) =
\{(\alpha_{2}, \beta_{2})\}$ and $B_{2}^{\beta_{3}}(C_{3},C_{6}) =
\{(\alpha_{2}, \beta_{3})\}$. Now suppose that $(\alpha_{0},
\beta_{0}) \approx_{2} (\alpha_{1}, \beta_{0})$. From \equref{p1'}
we obtain that $g \rhd_{0} h = g \rhd_{1} h$ for all $g \in
C_{3}$, $h \in C_{6}$ which is a contradiction. Thus
$B_{2}^{\beta_{0}}(C_{3},C_{6}) = \{(\alpha_{0}, \beta_{0}),
(\alpha_{1}, \beta_{0})\}$.
\end{remark}

\end{document}